\documentclass[a4paper,12pt]{amsart}
\usepackage{amssymb}
\usepackage{float}
\usepackage{graphicx}
\graphicspath{ {./images/} }
\restylefloat{figure,table}
\usepackage[colorlinks=true]{hyperref}
\hypersetup{allcolors=[rgb]{0.9,0.1,0.4}}
\usepackage[sort&compress,numbers]{natbib}
\usepackage{tikz}
\usetikzlibrary{positioning}
\usepackage{tikz-cd}
\usepackage{mathrsfs}
\usepackage{latexsym}
\usepackage{microtype}
\usetikzlibrary{decorations.markings,%
  matrix,%
  snakes}

\renewcommand{\H}{\mathcal{H}}

\theoremstyle{definition}
\newtheorem{definition}{Definition}

\newtheorem*{remark}{Remark}

\title{Higher Algebraic Structures and General Field Theories}
\author{Nils A.\ Baas}
\address{Department of Mathematical Sciences, NTNU, N-7491 Trondheim,
  Norway}
\email{nils.baas@ntnu.no}
\date{August 28, 2023}

\begin{document}

\begin{abstract}
    In this paper we show how the hyperstructure concept leads to new algebraic structures and general field theories.
\end{abstract}

\maketitle

\tableofcontents

\section{Introduction}
\label{sec:intro}

Higher structures appear in many ways in mathematics and theoretical and applied science. In a series of papers we have built up a framework for studying higher structures, see \cite{B1,B2,B3,B4,B5,B6,B7} and the references in there. We call the general notion of higher structures for Hyperstructures. Like in Topological Quantum Field theory there is often a need to represent Hyperstructures in other types of structures, for example algebraic ones, in order to assign states, observables and other properties. To find the right kind of algebraic structures is a subtle task.

Our goal in this paper is to introduce a new type of algebraic structures in which higher structures may be represented. Furthermore, this leads also to the suggestion of a type of General Field Theory. We will discuss these field theories and how they relate to ``local to global'' processes and organization of states. The main purpose of this paper is to sketch the main ideas.
\section{Higher Structures}
\label{sec:higher}

By higher structures we will mean hyperstructures in the sense we have defined and introduced them in \cite{B1,B2,B3,B4,B5,B6,B7}, where the $\H$-principle is given - the general basic principle in forming hyperstructures. Hyperstructures are designed to take care of situations with systems of systems of ..., structures of structures of..., etc.

Hyperstructures represent an extension of higher cobordisms and higher categories. The idea is that we start with a set or collection of objects $X$. We consider subsets of $X$, $S\subset X$ or $S \in \mathcal{P}(X)$, and 
$$ \Omega: \mathcal{P}(X) \longrightarrow \text{Set} $$ 
assigns a set (space) $\Omega(S)$ of properties for $S$.

Sometimes we may think of $\Omega$ as a functor (presheaf), but this is not always the case. $\Omega$ may be thought of partly as an observational mechanism of states, but may also include external fields or mechanisms acting on $S$. Hence, $\omega \in \Omega(S)$ may have both an internal and external part. We will consider pairs $(S,\omega)$ and form
$$ \Gamma = \{(S,\omega) \; \vert \; S\in \mathcal{P}(X) \text{ and } \omega \in\Omega(S)\}. $$
Intuitively we may consider ways of binding the elements of $S$. Define an assignment
$$ B: \Gamma \longrightarrow \text{Set} $$
\textemdash\ not necessarily a functor, but intuitively like $\Omega$ similar to a presheaf. We call $B(S,\omega)$ the set of bonds of $(S,\omega)$, extending the notion of morphism. $\text{Mor}(X,Y)$ extends to $B(X,Y,Z)$ etc. when more than two objects are involved. Cobordisms are bonds of the boundary components.

We iterate this construction by starting with
\begin{align*}
X_0,\; \Omega_0,\; \Gamma_0,\; B_0.
\end{align*}
The next level is then created by setting
\begin{align*}
X_1 = \{b_0 \; \vert \; b_0 \in B_0(S_0,\omega_0),\; S_0\in \mathcal{P}(X_0) \text{ and } \omega_0\in \Omega_0(S_0)\}.
\end{align*}

This is to be understood such that $S_0$ may vary and $X_1$ consists of elements from all non-trivial $B_0(S_0, \omega_0)$. Sometimes we reduce the set of interesting $S_0$'s due to some equivalence like $B_0(S_0, \omega_0) \simeq B_0(S_0', \omega_0')$. 

Depending on the situation we can now choose $\Omega_1$ and $B_1$ according to what we want to construct or study and then repeat the construction. This is not a recursive procedure since new properties and bonds arise at each level, and play a role in forming the next level.

It should be noted that according to the ``general principle'' in \cite{B7}, $X_1$ may be defined in many ways depending on various situations. For example like: $\prod_{S_0} B_0(S_0,\omega_0)$ or $\bigcup_{S_0}B_0(S_0, \omega_0)$ (disjoint union) or tensor product if applicable, or other formuli formed by $\{B_0(S_0, \omega_0)\}_{S_0}$. There are also situations where the bond spaces (sets) themselves are taken as bonds. Then we choose 
\[X_1 = \{B_0(S_0, \omega_0)\}_{S_0}.\]
With this in mind, the process continues in the usual way forming higher bonds and levels. 

Bonds extend morphisms in categories which are ``pairwise'' bonds. For graphs bonds are similar to ``edges'' in hypergraphs. Hence, a hyperstructure of order $n$ will consist of:

\begin{align*}
    \H: \begin{cases}
    X_0,\; \Omega_0,\; \Gamma_0,\; B_0 \\
    X_1,\; \Omega_1,\; \Gamma_1,\; B_1 \\
    \qquad \quad \vdots \\
    X_n,\; \Omega_n,\; \Gamma_n,\; B_n.
    \end{cases}
\end{align*}
At the technical level we require that
\begin{align*}
    &B_i(S_i,\omega_i) \cap B_i(S_i^\prime,\omega_i^\prime) = \emptyset
\end{align*}
for $S_i \neq S_i^\prime$ (``a bond knows what it binds'')
in order to define the $\partial_i$'s below, or we could just require the $\partial_i$'s to exist.

The levels are connected by ``boundary maps'' as follows:
\begin{align*}
    \partial_i: X_{i+1} \longrightarrow \mathcal{P}(X_i) \\
    \partial(b_{i+1}) = (S_i,\omega_i)
\end{align*}
and maps
\begin{align*}
    I_i: X_i \longrightarrow X_{i+1}
\end{align*}
such that $\partial_i \circ I_i = \text{id}$, giving $I_i$ as a kind of identity bond. 

Let us point out that the assignments by the $\Omega$'s and $B$'s may not be restricted to sets, but also spaces, manifolds, simplicial complexes, and algebraic structures in general or categories. Furthermore, $\mathcal{P}(X)$ may be replaced by more general constructions from $X$ – $Q(X)$ – where $Q(X)$ could be products, joins (if $X$ is a space), etc. 

For a further discussion of hyperstructures and examples we refer to \cite{B1,B2,B3,B4,B5,B6,B7}. Let us point out that in assigning or defining bonds of $(S_{i+1},\omega_{i+1}) = \{(S_i,\omega_i)\}$ the elements of $B_{i+1}(S_{i+1},\omega_{i+1})$ may depend on the $S_i$'s and $\omega_i$'s separately or in combination. ``Particles may be bonded through space, but also through states.''

\section{On multiplicative structures}
\label{sec:mstruct}

When we multiply we create complex objects from simpler ones. This is in a general way reflected in:

\begin{enumerate}
    \item[]
    \item $E_n$-algebras
    \item[]
    \item Operads
    \item[]
    \item Factorization algebras
    \item[]
    \item Higher TQFTs
    \item[]
    \item Haag-Kastler axioms
    \item[]
\end{enumerate}
We will sketch how they all fit into the general hyperstructure concept in a unifying way.

\subsubsection*{\underline{1) $E_n$-algebras}}
\text{ }\\
\text{ }\\
Let us follow J. Lurie (unpublished). $\mathcal{C}$ is a symmetric monoidal $\infty$-category, $n\geq0$. An $E_n$-algebra object of $\mathcal{C}$ consists of 
\begin{enumerate}
    \item[i)] $\forall \;U^{\text{open}}\! \subset \mathbb{R}^n \; (U\!\approx\!($homeomorphism) an open disk$)$, assign an object $A(U)\in \mathcal{C}$
    \item[ii)] $U_n,\ldots,U_m$ disjoint disks in a large disk $V$, a multiplication map:
    \begin{align*}
        \mu: \bigotimes_{1\leq i\leq n} A(U_i) \longrightarrow A(V)
    \end{align*}
    $\mu$ compatible with composition, equivalence for $m=1$. 
\end{enumerate}
$E_n$-algebras in $\mathcal{C}$ form another $\infty$-category
\begin{align*}
    \text{Alg}_{E_n}(\mathcal{C}) \cong \text{Alg}_{E_1}(\text{Alg}_{E_1}(\ldots (\mathcal{C})\ldots),
\end{align*}
see \cite{11}.

\subsubsection*{\underline{2) Operads}}
\text{ }\\
\text{ }\\
Family of objects 
\begin{align*}
    \mathcal{O} = \{\mathcal{O}(n)\}
\end{align*}
basically such that we have a composition (or ``multiplication'' operation)
\begin{align*}
    \mathcal{O}(n) \times \mathcal{O}(k_1)\times \cdots \times \mathcal{O}(k_n) \longrightarrow \mathcal{O}(k_1+\cdots +k_n)
\end{align*}
after coming from algebras and operations
\begin{align*}
    A \otimes \cdots \otimes A \longrightarrow A.
\end{align*}

Bonds in a hyperstructure with operations as described in \cite{17} extends operads and their algebras. The collection of objects that the bonds bind represent the algebras. The semantics being in the objects, the syntax in the bonds.

\subsubsection*{\underline{3) Factorization algebras}}
\text{ }\\
\text{ }\\
$X$ top. space, $\mathcal{O}(X)$ the category of open sets, $\mathcal{C}$ a symmetric monoidal category with weak equivalence. Prefactorization algebra $\mathcal{F}$ on $X$ is given by
\begin{align*}
    \mathcal{F}: \mathcal{O}(X) \longrightarrow \mathcal{C}
\end{align*}
such that if
\begin{align*}
    U_1 \;\cup\; &U_2 \cup \ldots \cup U_n \subset V \\
    (&\text{disjoint})
\end{align*}
we have a morphism
\begin{align*}
    \mathcal{F}(U_1) \otimes \cdots \otimes \mathcal{F}(U_n) \longrightarrow \mathcal{F}(V)
\end{align*}
and if 
\begin{align*}
    U_1  \;&\cup \ldots \cup \; U_{n_i} \subset V_i \quad \text{and} \\
    V_i \; &\cup \ldots \cup\; V_k \subset W
\end{align*}
we have
\begin{equation*}
\begin{tikzcd}[column sep=small]
    \bigotimes\limits_{i=1}^k \;\bigotimes\limits_{j=1}^{n_i} \mathcal{F}(U_j) \arrow[rr] \arrow[rd] & & \bigotimes\limits_{i=1}^k \mathcal{F}(V_i) \arrow[ld] \\
    & \mathcal{F}(W). & 
\end{tikzcd}
\end{equation*}
If $\mathcal{F}$ satisfies an additional gluing condition (equalizer condition), it is called a \textit{Factorization algebra}.

\subsubsection*{\underline{4) Higher TQFTs}}
\text{ }\\
\begin{align*}
    &\text{Bord}_n = n\text{-dim bordism }(n,\infty) \text{ category} \\
    &\mathcal{C}_n = \text{an } (n,\infty) \text{ category of suitable type}
\end{align*}
A TQFT is then given by an $n$-functor
\begin{align*}
    Z\!: \; \text{Bord}_n \longrightarrow \mathcal{C}_n.
\end{align*}
If $W^i$ is an $i$-dim manifold with boundary
\begin{align*}
    \partial W^i = \underbrace{W_1^{i-1} \cup \ldots \cup}_\text{in} \underbrace{W_{k_j}^{i-1,\ast}}_\text{out}
\end{align*}
\begin{align*}
    Z(W^i) \in Z(\partial W^i)=Z(W_i^{i-1}) \otimes \cdots \otimes Z(W_{k_{j-1}}^{i-1})^\ast
\end{align*}
indicating an upward type of ``multiplication'' in $Z$. This is what a higher TQFT should do.

\subsubsection*{\underline{5) Haag-Kastler axioms}}
\text{ }\\
\text{ }\\
In the Haag-Kastler axioms in mathematical physics one assigns to a domain in space-time $U$ a $C^\ast$-algebra $A(U)$. A nested structure of the domain may naturally give rise to a ``multiplicative'' structure of the $A(U)$-algebras similar to factorization algebras.

Let us add that there is also a concept of ``algebraic hyperstructures'' which is simply a set $X$ with an operation
\begin{align*}
    \ast: X \times X \longrightarrow \mathcal{P}(X)
\end{align*}
$\mathcal{P}(X)$ \textemdash\ the power set of $X$, non-empty and non-empty subsets. This is obviously covered by our notion of a $1$-level hyperstructure where
\begin{align*}
    B((x,y)) = x\ast y
\end{align*}
for ordered pairs $(x,y)$, otherwise the bond set is empty (see Wikipedia: Hyperstructures). Furthermore, Ran spaces, Configuration spaces and Clique complexes are structures that all can be suitably extended to higher versions.

All these structures fit into the conceptual framework of hyperstructures that we have introduced and studied in a series of papers, \cite{B1,B2,B3,B4,B5,B6,B7}. We will here suggest how to unify them in hyperstructures, which also paves the way for higher versions. Ultimately, this will lead us to the definition of a general field theory.

Starting with Operads and Props the ``Hyperstructure principle'' leads naturally to versions of Higher Operads and Higher Props. 

First, we will apply the ``Hyperstructure Principle'' to introduce higher states and observables in quantum field theory. Let $X_1$ be a given space (topological, manifold, stratified, space-time,...)

\begin{enumerate}
    \item[ ] $\mathcal{A}_1 = \{A_1(U_1)\; \vert\; U_1$ open in $X_1$ and $A_1$ an assignment of local states or observables, typically an algebraic assignment $\}$
    \item[ ] $X_2 =$ a space similar to $X_1$ (possibly $=X_1$)
\end{enumerate}
as for $\mathcal{A}_1$ we put:
\begin{enumerate}
    \item[ ] $\mathcal{A}_2= \{A_2(U_2)\; \vert\; U_2$ open in $X_2$, $A_2$ an assignment like $A_1\}$
    \item[ ] $\mathcal{P}(\mathcal{A}_1) = \{\text{(finite) subsets of objects}\}$
\end{enumerate}

For example,
\begin{align*}
    &S_1 \in \mathcal{P}(\mathcal{A}_1) \text{, then} \\
    &S_1 = \{A_1(U_1^1), \ldots ,A_1(U_1^{k_1}))\}
\end{align*}

We generalize multiplicative structures through ``bond'' formation:
\begin{align*}
    &B_1: \mathcal{P}(\mathcal{A}_1) \longrightarrow \mathcal{A}_2 \\
    &B_1(S_1) = B_1(\{A_1(U_1^1, \ldots , A_1(U_1^{k_1})\}) \in \mathcal{A}_2
\end{align*}
Hence, $B_1(S_1) = A_2(U_2)$. We continue this construction with $X_3, \mathcal{A}_3$ and
\begin{align*}
    &B_2: \mathcal{P}(\mathcal{A}_2) \longrightarrow \mathcal{A}_3 \\
    &A_3(U_3) \in B_2(\{A_2(U_2^1, \ldots , A_2(U_2^{k_2})\}) \\
    &U_3 \in X_3 \\
    &\quad\vdots \\
    &B_n: \mathcal{P}(\mathcal{A}_n) \longrightarrow \mathcal{A}_{n+1},\quad U_{n+1}\in X_{n+1}.
\end{align*}
Often $X_1=X_2=\ldots =X_{n+1}$ and sometimes $A_1=A_2=\ldots=A_{n+1}$, but we find it useful to formulate the general scheme.

Elements in the $A$ assignments we think of as:
\begin{align*}
    &\text{1st order states/observables } A_1(U_1),\quad U_1\subset X_1 \\
    &\text{2nd order states/observables } A_2(U_2),\quad U_2\subset X_2 \\
    &\quad\vdots \\
    &\text{$n$th order states/observables } A_n(U_n),\quad U_n\subset X_n \\
\end{align*}
often we may take $X_1=X_2=\ldots =X_n$ and consider nested open sets. In such a case the assiciated tensor products may be considered as bonds. The higher states/observables bind the lower ones. However, the point we are making is that ``multiplication'' (factorization $\otimes$-prod) brings us from one level to a higher. ``A product is a 'bond' of its factors''. This applies to operations in general. 

The bonds connect the $A_i$'s. In general hyperstructures going from the lower/local level to the higher/global are represented by what we call a ``globalizer'' \textemdash\ a generalized ``multiplicative'' structure which we will return to, see \citet{B6}.

Hyperstructures is a very general concept motivated by higher cobordisms and manifolds with singularities, and extending higher categories \textemdash\ higher morphisms being replaced by higher bonds. In addition the Obs assignment is new, see \cite{B1,B2,B3,B4,B5,B6,B7} and the references therein. The point we are trying to make here is that hyperstructures may be viewed as encompassing general multiplicative structures in a unifying way. For example unifying the structures in the Heisenberg and Schrödinger pictures in quantum theory.

If $\mathcal{H}$ and $\mathscr{S}$ are two hyperstructures we call an assignment
\begin{align*}
    \mathcal{F}: \mathcal{H} \longrightarrow \mathscr{S}
\end{align*}
a \textit{General Field Theory}. Intuitively we think of $\H$ representing space or system and $\mathscr{S}$ states or observables.

For applications outside mathematics and physics we keep the treatment here intentionally at a general level. We could tighten up the definitions making hyperstructures into some kind of higher categories and consider functors between them, see \citet{B7}. We postpone this for the time being.

In shorthand notation we have
\begin{align*}
    \H &= \{B_0,B_1, \ldots ,B_n  \} \\
    \mathscr{S} &= \{\mathscr{S}_0,\mathscr{S}, \ldots, \mathscr{S}_n  \}
\end{align*}
suppressing the rest of the notation. Consider a bond at level $i$,
\begin{align*}
    b_i \in B_i,
\end{align*}
then we ``dissolve'' it by the ``boundary'' map
\begin{align*}
    \partial_i: B_i &\longrightarrow B_{i-1} \\
    \partial_i(b_i) &= \{b_{i-1}^{k_{i-1}}  \}
\end{align*}
which leads us to think of $b_i$ being ``multiplicatively'' produced by the family $\{b_{i-1}^{k_{i-1}}\}$. But from a field theoretic point of view it may be more natural to consider this in the assignment hyperstructure $\mathscr{S}$ where we then think of $\mathscr{S}_n$ as the local level and $\mathscr{S}_0$ the top level.

In \citet{B3,B6} we described this as follows:
\begin{align*}
    \mathcal{F}_i: B_i \rightsquigarrow \mathscr{S}_{n-i}
\end{align*}
Global bonds are ``covered'' as follows:

\begin{equation*}
\begin{tikzcd}
    \{b(i_n)\} \arrow[r,"\partial_{n-1}"] & \{b(i_{n-1},i_n)\} \arrow[r, "\partial_{n-2}"] & \cdots \arrow[r,"\partial_0"] & \{b(i_0,\ldots,i_n)\}.
\end{tikzcd}
\end{equation*}
Clearly we want to connect the level assignments
\begin{align*}
    \mathcal{F}_k(\{b(i_k,\ldots ,i_n)   \})
\end{align*}
which we consider as a kind of ``level presheaves''. In order to do so we introduced in \citet{B3} a Groethendieck-type topology for the whole hyperstructure in order to glue the local level assignments together to a global object in $\mathscr{S}_0$. Then we get level connections $\delta_i$ dual to $\partial_{n-i}$ such that

\begin{equation*}
\begin{tikzcd}
    \mathcal{F}_n(\{b(i_n)\})  & \mathcal{F}_{n-1}(\{b(i_{n-1},i_n)\}) \arrow[l, "\delta_1"']  & \cdots \arrow[l, "\delta_2"']  & \mathcal{F}_0(\{b(i_0,\ldots,i_n)\}). \arrow[l, "\delta_n"']
\end{tikzcd}
\end{equation*}
The $\delta_i$'s may in general be ``relational'', but if we assume that we have suitable gluing conditions (similar to (pre)-sheaves) the $\delta_i$'s are functional and tell us how to get global states/observables from local assignments. We glue within and across levels, again see \cite{B3,B6}.

$\mathcal{F}$ is similar to a presheaf on $\H$ and when all $\{\delta_i \}$ exist we get a unique global state/observable from local ones making $\mathcal{F}$ into a ``sheaf-like'' assignment.

We call $\Delta=\{\delta_i\}$ a globalizer of the field theory
\begin{align*}
    \mathcal{F}: \H \rightsquigarrow \mathscr{S}.
\end{align*}
This makes our point that $\Delta$ may be viewed as a multilevel multiplicative structure
\begin{align*}
    \delta_{n-k}: \mathcal{F}_k(\{b(i_k, \ldots, i_n)\}) \longrightarrow \mathcal{F}_{k+1}(\{b(i_{k+1}, \ldots, i_n)\}).
\end{align*}
For if
\begin{align*}
    \partial_i(b_i) = \{b_{i-1}^{k_{i-1}}\}
\end{align*}
the globalizer will give an assignment

\begin{equation*}
\begin{tikzcd}
    \prod \limits_{k_{i-1}} \mathcal{F}_{i-1}(b_{i-1}^{k_{i-1}}) \arrow[r, "\delta_{n-i+1}"] & \mathcal{F}_i(b_i).
\end{tikzcd}
\end{equation*}

If there is a suitable tensor product it should be required that 
\begin{align*}
    \bigotimes_{k_{i-1}} \mathcal{F}_{i-1}(b_{i-1}^{k_{i-1}}) \longrightarrow \mathcal{F}_i(b_i).
\end{align*}
\section{General Field Theory}
\label{sec:GFT}

Let us be given an object or system we want to investigate or act on. Often it may come with a natural hyperstructure or we may impose one. We let the object or system be reperesented by a hyperstructure $\H$. For investigation or action it is then useful to assign some states, observables or other properties. We collect these in another hyperstructure $\mathscr{S}$. In general $\H$ and $\mathscr{S}$ may be arbitrary hyperstructures.

Assignments of the form
\begin{align*}
    \mathcal{F}: \H \rightsquigarrow \mathscr{S}
\end{align*}
is then what we will call a \textit{General Field Theory (GFT)}. This covers quantum field theories and may expand and enrich their structures. Also biological and social structures may be described by GFTs.

Most systems are built up of basic units that we consider elementary or local. Often we would like to deduce the global states and properties of the system from local ones. This may capture certain global states and properties, but not all. We suggest that one reason may be that such a direct deductional procedure ignores intermediate organizations into levels \textemdash\ namely what we have introduced as a Hyperstructure on the system of basic units. The level organization creates many new states and properties in all kinds of systems \textemdash\ physical as well as abstract ones. 

Often we want to act on a system in order to change its states or properties. In order to do so we may use another field theory
\begin{align*}
    \mathcal{G}: \H \rightsquigarrow \mathcal{E}
\end{align*}
where $\mathcal{E}$ acts on $\mathscr{S}$:
\begin{align*}
    \mathcal{E} \times \mathscr{S} \longrightarrow \mathscr{S}.
\end{align*}
$\H$ is often ``geometric'', $\mathscr{S}$ and $\mathcal{E}$ more ``algebraic'' in order to facilitate actions. High algebraic complexity may make more and interesting actions possible. Therefore higher algebraic hyperstructures as we have introduced are interesting and useful. 

In \cite[Section 5.2]{B5} we have considered hyperstructures of states $\mathscr{S}$ and pointed out the importance of finding a hyperstructure $\mathcal{E}$ acting on $\mathscr{S}$, 
$$\mathcal{E}\times \mathscr{S}\longrightarrow \mathscr{S},$$ 
such that the action may change global states through levelwise changes. This is important both in physical and biological systems. In neuroscience this corresponds to multilevel decoding acting on neurons to change global cognitive states. Similarily actions on genes to achieve organismic changes and actions on subatomic (subnuclear) levels in physical systems to achieve global atomic changes. 

Depending on the choice of $\mathscr{S}$ we get algebraic, geometric, topological, physical (quantum), biological (for example neural, genomic), etc. field theories. $\H$ represents the objects of study reflecting their basic nature. Then a GFT relates the two in an assignment
\begin{align*}
    \mathcal{F}: \H \rightsquigarrow \mathscr{S}.
\end{align*}

Assignments like
\begin{align*}
    \mathcal{F}: \H \rightsquigarrow \mathscr{S}
\end{align*}
may be conceived as a bond between the hyperstructures $\H$ and $\mathscr{S}$. This observation suggests a GFT may be taken to be a bond of families of hyperstructures $\{\H_i\}$
\begin{align*}
    B(\{\H_i\})
\end{align*}
which again may be organized into a hyperstructure of hyperstructures
\begin{align*}
    \text{Hyp}(B_j\{\H_{ij}\}).
\end{align*}
This is the ultimate General Field Theory, see also \citet{B5,B6}, emerging in stages
\begin{equation*}
\begin{tikzcd}
    \fbox{$\mathcal{F}: \H \longrightarrow \mathscr{S}$}\quad \arrow[d, bend right=50, shift right=11] \\
    \fbox{$B(\{\H_i\})$}\qquad\quad \arrow[d, bend right=50, shift right=10] \\
    \fbox{$\text{Hyp}(B_j(\{\H_{ij}\}))$}
\end{tikzcd}
\end{equation*}
representing a unification procedure of the $\H_i$'s.

For example in the family $\H_i$, one particular $\H$-structure may be a spatial collection of particles. Then one may design a whole family $\{\H_i\} \ni \H_i$ such that $B(\{\H_i\})$ (and possibly more levels) performs in a certain way \textemdash\ for example fusing particles (even nuclear!). Then we would call $B$ a \textit{fusion field} \textemdash\ a bond of a family of hyperstructures. Again, $B$ may be built into a total hyperstructure Hyp. This allows for example for more general interactions of objects and states.

In conclusion, GFTs are important in creating systems with new states and properties \textemdash\ in addition to analyzing existing systems. Furthermore, a field theory like 
$$\mathcal{F}\colon \H \rightsquigarrow \mathscr{S}$$
is useful in guiding how changes in $\H$ make changes in $\mathscr{S}$ and vice versa, like coding and decoding in neuroscience, see \cite{16}. 
See also \cite{B2,B3,B5,B6}.

\section{Multimodules}
\label{sec:mmod}

A topological quantum field theory is basically a functor
\begin{align*}
    Z: \mathcal{B} \longrightarrow \text{Vect}
\end{align*}
where $\mathcal{B}$ is a suitable cobordism category and Vect a suitable linear category of vector spaces (finite or infinite dimensional) representing states. Factorization algebras are roughly functors
\begin{align*}
    \mathcal{F}: \mathcal{O}(X) \longrightarrow \text{Alg}
\end{align*}
where $\mathcal{O}(X)$ is the category of open sets in a space $X$ and Alg a suitable category of algebras (observables). The first case is often called the Schrödinger picture while the second is called Heisenberg picture.

However, one wants to extend TGFTs to higher cobordism categories, and then we need interesting algebraic higher versions of Vect. The same idea applies to factorization algebras if we extend $\mathcal{O}(X)$ to ``higher order'' nested families of open sets. Higher versions of potential recipient categories have been constructed, see \cite{11,12,18}.

Higher cobordisms are really examples of hyperstructures binding boundary components, etc. Similarly are nested open sets examples of bonds in a hyperstructure. This is our motivation to extend TQFTs and Factorization algebras to the setting of hyperstructures. We think of a given hyperstructre $\H$ \textemdash\ for example representing the geometry or topology of the situation like cobordisms and nested open sets.

Then we want to assign (not necessarily functorial) states, observables and properties to the objects in $\H$. This assignment will need a similar recipient hyperstructure $\mathscr{S}$ \textemdash\ extending present versions of $\text{Alg}_n(S)$, $S$ some higher algebraic category, see \cite{11,12,18}.

Our goal is now to present some examples of such hyperstructures. In the two dimensional case categories of bimodules have been studied, and extended to higher categories of $\text{Bi}(\text{Bi}(\text{Bi}\ldots )))$-modules, see \cite{12}.

A bi-module is an abelian group $M$ with two actions \textemdash\ left and right \textemdash\ of two rings or algebras \textemdash\ $R$ and $S$, extended to derivators in categories, but the ideas are the same. In order to have a hyperstructure of this algebraic type representing other hyperstructures we need to extend this picture to tri-modules, four-modules and $n$-modules as a first step and then to have general, parametrized families of rings or algebras acting on $M$.

Preliminary thoughts: Candidate structure $\mathscr{S}=\mathcal{C}=\{\mathcal{C}_0,\ldots,\mathcal{C}_n\}$

\subsubsection*{I) \underline{Bimodules}}
\begin{align*}
    &\mathcal{C}_0 = \text{ rings or algebras}\\
    &\mathcal{C}_1 = \text{ abelian groups on which we have a pairwise interaction} \\
    & \qquad \qquad (R_1,R_2) \times M \longrightarrow R_1 \otimes M \otimes R_2 \approx M \\
    &\qquad \qquad  \qquad = \text{ bimodules} \\
    &\mathcal{C}_2 = \text{ bimodules of bimodules} \\
    &\qquad \vdots \\
    &\mathcal{C}_n \\
    &\text{Diagrammatic:} \\
    & \\
    & \qquad\qquad\qquad \begin{tikzpicture}
\node[circle, minimum size=4mm, draw=black](1){$1$};
\node[circle, minimum size=12mm, draw=black](e)[right=of 1]{};
\node[circle, minimum size=4mm, draw=black](2)[right=of e]{$2$};
\draw (1) -- (e);
\draw (e) -- (2);
\end{tikzpicture}
\end{align*}

\subsubsection*{II) \underline{Tri-modules}}
\begin{align*}
    &\mathcal{C}_0 = \text{ as before} \\
    &\mathcal{C}_1 = \text{ triple interaction} \\
    & \qquad \qquad (R_1, R_2, R_3) \times M \longrightarrow M \\
    & \qquad \qquad \text{Various conditions may come in here, like the action being } \\
    & \qquad \qquad \text{pairwise commutative, or more refined conditions.} \\
    & \qquad \qquad \qquad =\text{tri-modules} \\
    & \mathcal{C}_2 = \text{ tri-modules of tri-modules...} \\
    &  \vdots \\
    & \mathcal{C}_n \\
    & \text{Diagrammatic:} \\
    &\qquad\qquad\qquad \begin{tikzpicture}
\node[circle, minimum size=4mm, draw=black,](1){$1$};
\node[circle, minimum size=12mm, draw=black,](e)[right=of 1]{};
\node[circle, minimum size=4mm, draw=black,](2)[above=of e]{$2$};
\node[circle, minimum size=4mm, draw=black,](3)[right=of e]{$3$};
\draw (1)--(e)--(3);
\draw (e)--(2);
\end{tikzpicture}
\end{align*}

\subsubsection*{III) \underline{$m$-modules}}
\begin{align*}
    & \mathcal{C}_0 = \text{ as before} \\
    & \mathcal{C}_1 = m \text{ element interaction} \\
    & \qquad \qquad (R_1,R_2, \ldots, R_m) \times M \longrightarrow M, \\
    & \qquad \qquad \text{for example pairwise commutativity,or a suitable action of} \\
    & \qquad \qquad \text{$\mathcal{R}=\bigotimes_{i=1}^m R_i,$ or actions connected through a parameter space.} \\
    &  \vdots \\
    & \mathcal{C}_n \\
    & \text{Diagrammatic:} \\
    & \qquad\qquad\qquad \begin{tikzpicture}
\node[circle, draw=black, minimum size=4mm](1)at(-0.8,1.7){$1$};
\node[circle, draw=black, minimum size=4mm](2)at(0,2){$2$};
\node[circle, draw=black, minimum size=4mm](3)at(0.8,1.7){$3$};
\node[circle, draw=black, minimum size=12mm](e)at(0,0){};
\node[circle, draw=black, minimum size=4mm](m)at(2,0.1){$m$};
\draw (1)--(e)--(2);
\draw (3)--(e)--(m);
\draw[dotted, thick] (1.2,1.3)--(1.7,0.6);
\draw[dotted, thick] (0.7,0.7)--(1,0.2);
\end{tikzpicture}
\end{align*}

By iteration we get $m$-modules with $n$ levels. The cardinality of the acting families could be the same or vary from level to level (depending on Obs in the design). Let $m=(m_1,m_2,\ldots,m_k)$. Then we form
\begin{align*}
    \mathcal{C}(m_1,m_2,\ldots,m_k)
\end{align*}
forming $m_{i+1}$-modules from $m_i$-modules, for example by suitable tensor product actions. 

This is the motivating picture which we will now extend to more general families of algebraic objects acting on $M \in \mathcal{M}$ \textemdash\ some chosen basic structures, for example abelian groups. As an example we may consider multidimensional matrices (hypermatrices), see \cite{19}, and higher levels created by matrices of such matrices etc.

An example of a multimodule structure: Let 
\begin{center}
\begin{tikzcd}
    M_{i_0}^1 \arrow[r, shift left=3] \arrow[r, "\vdots", shift right=2] & M_{i_0 i_1}^2 \arrow[r, "\vdots", shift right=2] \arrow[r, shift left=3] & \cdots \arrow[r, "\vdots", shift right=2] \arrow[r, shift left=3] & M_{i_0 \cdots i_n}^n
\end{tikzcd}
\end{center}
be a reverse tree of smooth maps and manifolds. The rings (algebras) of $C^\infty$ fuctions are then connected as 
\begin{center}
    \begin{tikzcd}
        C^\infty (M_{i_0}^1) & C^\infty(M_{i_0 i_1}^2) \arrow[l, shift right=3] \arrow[l, "\vdots"', shift left=2] & \cdots \arrow[l, "\vdots"', shift left=2] \arrow[l, shift right=3] & C^\infty(M_{i_0 \cdots i_n}^n) \arrow[l, "\vdots"', shift left=2] \arrow[l, shift right=3]                 
        \end{tikzcd}
\end{center}
Let $w=(i_0\cdots i_{k-1})$ and $R(w) = C^\infty(M_{i_0 \cdots i_k}^k).$ Then $R(w,i_k)\longrightarrow R(w)$ and we have actions 
$$\bigotimes_{i_k}R(w,i_k)\times R(w)\longrightarrow R(w).$$

Let us start with a family of ``structures'' \textemdash\ for example rings or algebras.
\begin{align*}
    \mathscr{S}_0 = \{S_{t_0} \;\vert \; t_0 \in T_0, T_0 \text{ a space, for example a manifold} \}
\end{align*}
Given another structure $S_0\in \mathcal{M}_0$ ($\mathcal{M}_0$ as $\mathcal{M}$ above).
\begin{definition}
$S_0$ is called a module over $\mathscr{S}_0$ iff we have an action
\begin{align*}
    \mathscr{S}_0 \times S_0 \longrightarrow S_0
\end{align*}
satisfying reasonable conditions. $\mathscr{S}_0$ may be finite, infinite or uncountable.
\end{definition}

For example: $\mathscr{S}_0 = \{R_{t_0} \;\vert \; t_0 \in T_0\}$, where $R_{t_0}$ is a ring or algebra or another algebraic structure, and $T_0$ is a parameter space, for example a manifold.

Next we consider families of actions parametrized by a space (manifold) $W_0$:
\begin{align*}
    W_0 \times \mathscr{S}_0 \times S_0 \longrightarrow S_0.
\end{align*}
Then we proceed to the next level:

\noindent Let $\mathscr{S}_1 =$ \{family of modules over $\mathscr{S}_0$, parametrized by $T_1$\}. 

\noindent Let $S_1\in \mathcal{M}_1$ -- a structure like $\mathcal{M}_0$.

\begin{definition}
$S_1 \in \mathscr{M}_1$ is called a module over $\mathscr{S}_1$ iff we have an action
\begin{align*}
    W_1 \times \mathscr{S}_1 \times S_1 \longrightarrow S_1,
\end{align*}
$W_1$ a space parametrizing the actions.
\end{definition}
Proceeding inductively:
\begin{align*}
    W_i \times \mathscr{S}_i \times S_i \longrightarrow S_i
\end{align*}

\begin{definition}
$\mathscr{S} = \{\mathscr{S}_0,\mathscr{S}_1, \ldots,\mathscr{S}_n\}$ is a hyperstructure of modules extending iterated bi-modules and $m$-modules.
\end{definition}

The ``boundary'' maps in the hyperstructure are given by
\begin{align*}
    \partial: \mathscr{S}_{i+1} &\longrightarrow \mathscr{S}_i \\
    \partial(S_{i+1}) &= \{S_i^{t_i} \;\vert \; t_i \in T_i, \text{ acting on } S_{i+1} \}.
\end{align*}
$S_{i+1}$ binds the family $\{S_i^{t_i} \}$ and Obs is detecting the actions.

\noindent What we have achieved so far is the following: \newline
\noindent Starting with some basic objects $\mathcal{M}_0$ (like abelian groups) and a family of structures (rings, algebras,...) $\mathscr{S}_0$, we have created a hyperstructure of level $n$
\begin{align*}
    \mathscr{S} = \{ \mathscr{S}_0,\mathscr{S}_1,\ldots,\mathscr{S}_n \}.
\end{align*}
We call $\mathscr{S}$ a multimodule structure (more general than a higher category).

Included in the construction of the actions are parameter spaces
\begin{align*}
    (T_0,W_0), (T_1,W_1), \ldots ,(T_n,W_n)
\end{align*}
which we may assume to be manifolds. This structure is a vast extension of iterated bi-module structures and their higher categories.

Let us consider higher cobordisms (up to a certain dimension) as bonds in a hyperstructure $\H_{\text{Bord}}$. Moreover, if $X$ is a topological space $\H_{\text{Nest}}$ is the hyperstructure of nested open sets. Bonds being inclusions, Obs checking openness. Then assignments
\begin{align*}
    \H_{\text{Bord}} \rightsquigarrow \mathscr{S}=\H_{\text{Alg}}
\end{align*}
and
\begin{align*}
    \H_{\text{Nest}} \rightsquigarrow \mathscr{S}=\H_{\text{Alg}}
\end{align*}
will generalize TQFTs, Factorization Algebras and Factorization Homology.

Let us now see how a stratified space leads to a multimodule structure. \newline
\noindent $Y =$ a stratified space or a manifold with a stratification, $\mathscr{S}_0,\mathcal{M}_0$ basic algebraic objects. The desired parameter spaces come from the stratification. $Y$ has a filtration such that
\begin{align*}
    Y_0 \subset Y_1 \subset \ldots \subset Y_i\subset \ldots \subset Y_n = Y
\end{align*}
and locally
\begin{align*}
    Y_i - Y_{i-1} \sim \mathbb{R}^i\times CT_i,
\end{align*}
$T_i$ of dimension $n-i-1$. The local structure may extend to a ``neighborhood'' of the singular set like
\begin{align*}
    W_i \times CT_i,
\end{align*}
$W_i$ a manifold of dim $i$. For example consider two $2$-spheres intersecting each other, a neighborhood of singular set $\approx S^1 \times C\mathbb{Z}_4$. $T_i$ may be a stratified space again so the process may be iterated as in Figure 1:

\begin{figure*}[hbt!]
\begin{equation*}
\begin{tikzpicture}
\draw (0,0) -- (1,1) -- (2,0);
\draw[dashed] (2,0) arc (0:180:1 and 0.1);
\draw (0,0) arc (180:360:1 and 0.1);
\draw (1,2.5) -- (2,3.5) -- (3,2.5);
\draw[dashed] (3,2.5) arc (0:180:1 and 0.1);
\draw (1,2.5) arc (180:360:1 and 0.1);
\draw (1.5,3.6) -- (2.5,4.6) -- (3.5,3.6);
\draw[dashed] (3.5,3.6) arc (0:180:1 and 0.1);
\draw (1.5,3.6) arc (180:360:1 and 0.1);
\draw[dashed] (1,1.1) -- (1.4,2.4);
\draw (-0.8,0.6) -- (0,0);
\draw (0,3.2) -- (1,2.5);
\draw (0.6,4.3) -- (1.5,3.6);
\draw (1.7,-0.2) arc (180:290:0.4);
\draw (2.7,2.3) arc (180:290:0.4);
\draw (3.2,3.4) arc (180:290:0.4);
\node[]at(-1,0.4){\begin{small}$W_0$\end{small}};
\node[]at(-0.3,2.9){\begin{small}$W_{n-1}$\end{small}};
\node[]at(0.5,4){\begin{small}$W_n$\end{small}};
\node[]at(2.5,-0.5){\begin{small}$T_0$\end{small}};
\node[]at(3.7,2){\begin{small}$T_{n-1}$\end{small}};
\node[]at(4,3.1){\begin{small}$T_n$\end{small}};
\end{tikzpicture}    
\end{equation*}
Figure 1
\end{figure*}

This picture is also what we use in manifolds with ``joinlike''-singularities and their cobordisms \textemdash\ see \citet{B8}. Compare with the structures 
\begin{align*}
    A = \{A(\omega)\}, \quad \partial_i A(\omega) = A(\omega) \times P_i, \quad i \not\in \omega
\end{align*}
as described in \citet{B8}.

The $\mathscr{S}$ structure at level $i$ will then be given by
\begin{align*}
   \underbrace{W_i \times T_i}_\text{Topological}
   \times \mathscr{S}_i \times S_i \longrightarrow S_i
\end{align*}
We require the maps to be locally constant with respect to open sets in the strata of the stratification. This shows how a stratified space together with basic Algebraic structures $\mathscr{S}_0, \mathcal{M}_0$ lead to a Hyperstructure $\H_\text{Alg}$ which can now be used as a recipient of GFTs extending TQFTs and Factorization algebras.

A stratified space $Z$ may be considered as a hyperstructure where the bonds are the strata \textemdash\ the zero dimensional ones at the bottom and the top singularity(ies) the top bond. Notation $\H_\text{Strat}(Z)$. From the stratification we also make a hyperstructure of nested open sets \textemdash\ $\H_\text{Nest}(Z)$. 

If $Z$ is stratified, we get an assignment
\begin{equation*}
\begin{tikzcd}
   \H_\text{Nest}(Z) \arrow[r, rightsquigarrow, "\mathcal{F}"] & \H_\text{Alg}(Z)
\end{tikzcd}    
\end{equation*}
meaning that open sets at level $i$ go to the $\mathscr{S}_i$ in $\mathscr{S}$ constructed from $Z$'s stratification. This is a vast extension of a factorization algebra. $\H_\text{Nest}(Z)$ is the hyperstructure of nested open sets as follows:
\begin{align*}
    \mathcal{U}=U(i_1,\ldots,i_n)\} \quad \text{where} \quad U(i_1,\ldots,i_k,\ldots,i_n) \subset U(i_1,\ldots,\hat{i}_k,\ldots,i_n)
\end{align*}
and the larger open sets are bonds of the smaller. Furthermore
\begin{align*}
    \partial_j U(i_1,\ldots,\hat{i}_j,\ldots,i_n) = \{  U(i_1,\ldots,i_j,\ldots,i_n)  \}.
\end{align*}
Inclusions are bonds and Obs checks openness.

Cobordism categories (like Bord) may also be considered as hyperstructures and given representations in some general hyperstructures $\H$, $\H_\text{Alg}$ or some $\H_\text{Alg}(Z)$ (for a specific stratified space $Z$). This extends TQFTs.

$\H_\text{Strat}(Z)$ and $\H_\text{Nest}(Z)$ may to give field theories be represented in $\H_\text{Alg}(Z)$ or $\H_\text{Alg}(W)$ ($W$ an arbitrary stratified space) or a quite general $\H$.  We require the assignments to be locally constant with respect to open sets in the strata of the stratification.

$\mathcal{F} \colon \H_\mathrm{Geom} \rightsquigarrow \H_\mathrm{Alg}$ extends the construction of factorization homology and may symbolically be written as 
\begin{align*}
    \mathcal{F} = \int_{\H_\mathrm{Geom}}\!\! \H_\mathrm{Alg}.
\end{align*}
Let us think of $\H_\mathrm{Geom}$ as being represented by bonds $\{B_0,B_1,\ldots,B_n\}$. Then $\mathcal{F}$ will perform assignments in $\H_\mathrm{Alg}$ as:
\begin{align*}
    B_0 &\rightsquigarrow \mathrm{Alg}(B_0) \quad \text{--} \quad \mathscr{S}_n\text{-module} \\
    B_1 &\rightsquigarrow \mathrm{Alg}(B_1) \quad \text{--} \quad \mathscr{S}_{n-1}\text{-module} \\
    &\quad\vdots \qquad \qquad \qquad \qquad \vdots\\
    B_n &\rightsquigarrow \mathrm{Alg}(B_n) \quad \text{--} \quad \mathscr{S}_0\text{-module} 
\end{align*}
-- the $\mathscr{S}_i$'s are of the type previously described.

Stratifications may enter in two ways:
\begin{itemize}
    \item[(i)] in the geometric object to be represented in $\H_\mathrm{Alg}$.
    \item[(ii)] in the construction of $\H_\mathrm{Alg}$ itself.
\end{itemize}

They are all examples of General Field Theories of the type
\begin{equation*}
\begin{tikzcd}
   \H \arrow[r, rightsquigarrow, "\mathcal{F}"] & \mathscr{S}.
\end{tikzcd}
\end{equation*}
Note that stratifications also play a role in:
\begin{enumerate}
    \item[i)] Giving a stratified space a hyperstructure
    \item[ii)] Organizing algebraic actions in a hyperstructure
    \item[iii)] Generalizing Factorization Algebras and TQFTs
\end{enumerate}

Let us also make a comment on levels of observables and states. Ayala \& Francis \cite{9} show that one can compute global observables from locals:
\begin{align*}
    \int_M \text{obs}(\mathbb{R}^n) \longrightarrow \text{obs}(M).
\end{align*}
This map will be an isomorphism (equivalence) for a perturbative theory, but not for a non-perturbative. However, this may come from ignorance of possible higher structures in $M$, at least levels add more states. Suppose
\begin{align*}
    M \rightsquigarrow \H(M) = \{B_0,B_1,\ldots,B_n   \}
\end{align*}
in terms of bonds. Then we may have a picture as follows:
\begin{equation*}
\begin{tikzcd}
  \int\limits_{\H} \text{obs}(B_0) \arrow[r] \arrow[d, rightsquigarrow] & \int_\H \text{obs}(B_1)\arrow[r] \arrow[d, rightsquigarrow] & \cdots \arrow[r] & \int_\H \text{obs}(B_n) =\!\!\! \begin{tabular}{c}
      \text{Global}   \\
    \text{observables}
  \end{tabular} \arrow[d, rightsquigarrow, shift right=15] \\
  \mathscr{S}_0 \arrow[r] & \mathscr{S}_1 \arrow[r] & \cdots \arrow[r] & \mathscr{S}_n = \text{Global states}
\end{tikzcd}
\end{equation*}
The horizontal arrows represent a ``globalizer'' in the sense of \citet{B3,B6}.

Let us sum up the ideas presented here. Given a space $X$, we may often obtain a geometric or topological decomposition or filtration 
$$X(w, i_{n+1})\subseteq X(w) = X(i_1, i_2, \ldots, i_n)\subseteq X.$$

The object $\H(X) = \{X(w)\}$ may be thought of as a natural hyperstructure on $X$: 
$$\partial(X(w)) = \{ X(i_1,\ldots, i_n, i_{n+1}) \},$$
similar to manifolds with cone-singularities, as in \cite{B8}. It is often useful to associate states or properties to $\H(X)$ in the form of a general field theory 
$$\mathcal{F}\colon \H(X)\rightsquigarrow \text{Alg}.$$
We think of $\text{Alg}$ as an algebraic hyperstructure with objects $A=\{A(w)\}$. The bonds for example being via a multimodule structure based on a basic family of rings or algebras. This means, by using the same indexing on both sides, 
\begin{center}
    \begin{tikzcd}
        X(w) \arrow[d, "\partial"] \arrow[r, snake] & A(w) \arrow[d, "\partial"] \\
        {\{X(w,i)\}} \arrow[r, snake]                & {\{A(w,i)\}}              
        \end{tikzcd}
\end{center}
and the bonding by 
$$\bigotimes_i A(w,i)\times A(w)\longrightarrow A(w).$$

The idea here is to extend the construction of $\text{Alg}_{E_n}(\mathcal{C})$, see \cite{11}, to bonds and hyperstructures. 

In constructing $E_n$-algbras we use $\mathbb{R}^n$ as a model space with a hyperplane stratification. Here we suggest using a stratified space $Y$ where $A(w)$ is a locally constant assignment of $Y(w)$ to $\mathcal{C}$ which could be a symmetric monoidal $\infty$-category or some other type of hyperstructure, getting $\text{Alg}(\mathcal{C})$. Anyway, we think of $\text{Alg}(\mathcal{C})$ as a bond extension of $\text{Alg}_{E_n}(\mathcal{C})$. 
\section{Multistructures}
\label{sec_multistruct}

Generalized Field Theories appear in many areas of mathematics and science as we have discussed in \citet{B3,B5,B6}, and we have found that hyperstructures is an adequate framework for studying such level structures. Naturally the states, observables and properties of such structures are of great interest. GFTs represent a method to create new states, observe existing level structures and create fields to act on them for desired changes as presented in \citet{B5}. We will now go deeper into the underlying mathematical structures.

Let $X$ be an object to be studied and acted upon. The procedure we suggest is:
\begin{align*}
    &X \rightsquigarrow \H(X), \quad \text{put a hyperstructure on } X \\
    &\H(X) \rightsquigarrow \mathscr{S}, \quad \text{assign states, observables...} \\
    &\H(X) \rightsquigarrow \mathcal{E}, \quad \text{introduce an action field for changes} \\
    &\mathcal{E} \times \mathscr{S} \rightsquigarrow \mathscr{S}.
\end{align*}

We start with basic objects $\mathscr{S}_0,\mathcal{E}_0$ \textemdash\ later to be the bottom level of hyperstructures $\mathscr{S}, \mathcal{E}$. For example $\mathcal{E}_0$ could be ``waves'' acting on objects in $\mathscr{S}_0$ and $\mathcal{E}$ a hyperstructure of such waves. $\mathscr{S}_0$ could be vector spaces, groups, algebras, manifolds and $\mathcal{E}_0$ objects acting on these, often of similar type.

\textit{Organizing Principle:} When we have decided what the basic objects should be, how do we organize them for effective use and actions? Our suggestion is that we build $\mathscr{S}_0$ and $\mathcal{E}_0$ into suitable hyperstructures as bottom levels.

Let us sketch the general scheme we have in mind. Let $\mathcal{C}_0$ be a family of structures covering (including) the structures for $\mathscr{S}_0$ and $\mathcal{E}_0$. There is a plethora of ways to extend $\mathscr{S}_0$ to a full hyperstructure
\begin{align*}
    \mathscr{S} = \{\mathscr{S}_0,\ldots,\mathscr{S}_n  \}
\end{align*}
and similarly for $\mathcal{E}_0$. Think of $\mathcal{C}_0 = \mathscr{S}_0$.

\textit{Multistructures} (generalizing multimodules) based on $\mathscr{S}_0$.
\begin{align*}
    &\mathscr{S}_0 = \text{ basic structures} \\
    &\mathscr{S}_1 = \text{ objects such that families from } \mathscr{S}_0 \text{ act on them:} \\
    & \qquad \qquad \text{Fam}(\mathscr{S}_0) \times \text{obj.}(\mathscr{S}_1) \longrightarrow \text{obj.}(\mathscr{S}_1).
\end{align*}
Such objects are considered as \textit{bonds} for the families \textemdash\ extending the concept of bi- and multimodules.

In forming a hyperstructure (in standard notation) Obs of Fam is the action. This process goes on:
\begin{align*}
    &\mathscr{S}_2 = \text{ objects such that} \\
    &\qquad \quad \text{Fam}(\mathscr{S}_1) \times \text{obj.}(\mathscr{S}_2) \longrightarrow \text{obj.}(\mathscr{S}_2) \\
    & \quad \vdots \\
    &\mathscr{S}_n = \text{ objects such that} \\
    &\qquad \quad \text{Fam}(\mathscr{S}_{n-1}) \times \text{obj.}(\mathscr{S}_n) \longrightarrow \text{obj.}(\mathscr{S}_n).
\end{align*}
In general then
\begin{align*}
    \mathcal{C}_0 = \mathscr{S}_0\; \rightsquigarrow \; \H = \{\mathscr{S}_0,\mathscr{S}_1,\ldots,\mathscr{S}_n   \} 
    = \{\mathcal{C}_0,\mathcal{C}_1,\ldots,\mathcal{C}_n\}.
\end{align*}
We think of $\mathscr{S}$ as internal fields and $\mathcal{E}$ as external fields. It may be useful to construct $\mathcal{E}$'s with a high degree of complexity in order to perform desired changes of the internal $\mathscr{S}$'s. It is also useful when they have the same type of compatible structures.

Multimodules as we have introduced them here are clearly examples of multistructures. Hyperstructures of multimodules or multistructures may be used for representation of any kind of hyperstructure, in particular geometric and topological ones like in general quantum field theories.
\section{General algebraic operations and Fusion Processes}
\label{sec:fusion}

The $\H$-principle is about how many things act together to make new things of higher order. As we have pointed out ``multiplication'' in TQFTs and factorization algebras may naturally be extended to create higher order objects and an $\H$-structure. Intuitively and in general, ordinary multiplication leads to more complex objects, but one may also say this about addition.

Let us consider (algebraic) operations in general for sets $A$ with some given structure (algebraic). We think of $\Box$ as a general operation like multiplication, addition, tensor, etc. \textemdash\ binary or multiple.
\begin{align*}
    \Box_1: A_{i_1}^1 \times \cdots \times A_{i_n}^1 \longrightarrow A_{i^2}^2
\end{align*}
Then we use Obs for detecting properties (states) and proceed:
\begin{align*}
    \Box_1(\{ a_{i_k}^1, \omega_1 \in \text{Obs}_1  \}) \in A_{i^2}^2.
\end{align*}
The next level \textemdash\ even within the starting level \textemdash\ is then 
\begin{equation*}
\begin{tikzcd}
    \Box_2(\{ a_{i_k}^2, \omega_2 \}) \in A_{i^3}^3 \arrow[d,rightsquigarrow, "\text{ continuing}"]\\
    \Box_n(\{ a_{i_k}^n, \omega_n \}) \in A_{i^{n+1}}^{n+1}.
\end{tikzcd}    
\end{equation*}
This defines a higher, hyperstructured set of operations. In \citet{B10} we have considered ``addition'' of genes with properties (features) and shown how $3+2$ has different properties than $5$, see Figure 2.

\begin{figure}[hbt!]
\begin{equation*}
\begin{tikzpicture}
\draw (0,0) ellipse (1 and 0.5);
\draw (1.7,0) ellipse (0.6 and 0.5);
\draw (0.7,0) ellipse (1.9 and 0.7);
\node[circle]at(0,0){\begin{Large}$\cdots$\end{Large}};
\node[circle]at(1.7,0){\begin{Large}$\cdot\cdot$\end{Large}};
\node[circle]at(3,0){$\neq$};
\draw (4.8,0) ellipse (1.4 and 0.7);
\node[circle]at(4.8,0){\begin{Large}$\cdot\cdot\cdot\cdot\cdot$\end{Large}};
\end{tikzpicture}
\end{equation*}
Figure 2
\end{figure}

The idea is to enrich the structures and computations by adding properties levelwise and using them in the next level construction. This indicates how one could take pure and entangled states in quantum theory (like GHZ-states) to higher versions.

We may also think of multiple inputs and outputs, including multiplication with inverse, and the same for addition. One output may be a desired release, like for example energy in nuclear fusion and fission. General fusion processes may have this feature: many inputs, several outputs where \textit{one} or \textit{more} are desired.

This applies also to biological systems \textemdash\ ``fusion of genes'' and neural assemblies. In social science \textemdash\ subpopulations with creation or release of properties, giving new resources useful for protecting the structure (``defence''). If such operations are naturally available, this is an interesting way to produce fusion process \textemdash\ even of higher order.
Fusion:

\begin{figure}[hbt!]
\begin{equation*}
\begin{tikzpicture}
\node[circle, draw=black](X)at(0,1){$X$};
\node[circle, draw=black](Y)at(0,0){$Y$};
\draw (-0.7,-0.7) rectangle (0.7,1.7);
\node[circle](arrow)at(2,0.5){$\longrightarrow$};
\node[circle](X2)at(4,0.85){$X$};
\node[circle](Y2)at(4,0.15){$Y$};
\draw (4,0.5) ellipse (0.5 and 1);
\draw (3.3,-0.7) rectangle (4.7,1.7);
\end{tikzpicture}    
\end{equation*}
Figure 3
\end{figure}

\textit{General Fusion}.

If we want to fuse two objects or collections, we may facilitate this in two basic ways (see discussion in \citet{B2,B3}).
\begin{enumerate}
    \item[I)] Put an $\H$-structure on the (two) \textit{fusing objects} $X$ and $Y$ \textemdash\ $\H(X), \H(Y)$. \newline
    \noindent $\Box =$ just putting them together and realize:
    \begin{align*}
        \H(X\Box Y) = \H(X) \Box_\H \H(Y).
    \end{align*}
    Example: $X$ and $Y$ stratified spaces glued levelwise along strata, give a new stratified space.
    
    \noindent Similar for genes in reproduction and chromosome formation. This is a direct process \textemdash\ fusion in space.
    \item[II)] Like in general field theory we may have state assignments:
    \begin{equation*}
    \begin{tikzcd}
       \mathcal{F}: 
       \begin{tabular}{c}
         $\H(X), \H(Y)$ \\
        $\text{ or} $\\
        $\H(X\Box Y)$
    \end{tabular} 
    \arrow[r, rightsquigarrow] & \mathscr{S}
    \end{tikzcd}
    \end{equation*}
$\mathscr{S} =$ hyperstructure of states.

\noindent $S \in \mathscr{S}$ may be an unfused state and $S^\prime \in \mathscr{S}$ a derived fused state. Then we will provide an action
\begin{align*}
    A: \mathscr{S} \longrightarrow \mathscr{S}
\end{align*}
such that $A(S) = S^\prime$. This may be done by some dynamic process, bringing relevant inputs or by an action of another hyperstructure $\mathcal{E}$:
\begin{align*}
    \mathcal{E} \times \mathscr{S} \longrightarrow \mathscr{S}.
\end{align*}
$\mathcal{E}$ could for example be one of the types we have introduced here. 
\end{enumerate}
The important point is that ``producable'' actions on the local level may propagate through the hyperstructure to give the desired global action \textemdash\ like in this case of fusion of two or more objects. In understanding hyperstructures ``sociological'' metaphors are often useful.

In conclusion, we have introduced hyperstructures, given examples to be used in the study of general field theories and presented a setting for general field theories and fusion.
\section{Bonds and states}
\label{sec:bonds}

Given
\begin{align*}
    \H = \{B_0,B_1,\ldots,B_n\}, \text{ a hyperstructure}
\end{align*}
\begin{equation*}
\begin{tikzcd}
    \H \arrow[r, rightsquigarrow, "\mathcal{F}"] & \mathscr{S}
\end{tikzcd}, \text{ a GFT}
\end{equation*}
New (higher) bonds may be formed as follows:
\begin{enumerate}
    \item[i)] $b$ ($b_i$'s) bonds of existing bonds or objects
    \begin{equation*}
         \begin{tikzpicture}
\node[circle, draw=black, minimum size=4mm](bi)at(0,0){$b_i$};
\node[circle, draw=black, minimum size=8mm](1)at(2,0){};
\node[circle, draw=black, minimum size=8mm](2)at(4,0){};
\draw (2,0) ellipse (3 and 1);
\node[circle](b)at(3.9,1){$b$};
\end{tikzpicture}
    \end{equation*}
   
    \item[]
    \item[ii)] \textit{ }\\
    \begin{figure}[hbt!]
    \begin{equation*}
        \begin{tikzpicture}
\draw (-2,0) ellipse (0.5 and 1);
\draw (0,0) ellipse (0.5 and 1);
\draw (2,0) ellipse (0.5 and 1);
\draw[fill=black] (-2,0.7) circle (0.1);
\draw[fill=black] (0,0.7) circle (0.1);
\draw[fill=black] (2,0.7) circle (0.1);
\draw[fill=black] (0,2) circle (0.1);
\draw (-2,0.7) -- (0,2) -- (0,0.7);
\draw (0,2) -- (2,0.7);
\draw[fill=black] (-2,0.2) circle (0.1);
\draw[fill=black] (0,0.2) circle (0.1);
\draw[fill=black] (2,0.2) circle (0.1);
\draw (-2,0.2) -- (0,0.2) -- (2,0.2);
\draw (-2.45,-0.35) -- (-1.55,-0.35);
\draw (-0.45,-0.35) -- (0.45,-0.35);
\draw (1.55,-0.35) -- (2.45,-0.35);
\node[circle](bi)at(-2,-0.6){$b_i$};
\node[]at(-2.25,0){\begin{small}$\mathscr{S}_i$\end{small}};
    \end{tikzpicture}
    \end{equation*}
    Figure 4
    \end{figure}

forming bonds of \textit{states}, binding both lower bonds and states.
\end{enumerate}

If we are given a basic set of objects $X$, and can associate to it a hyperstructure of states $\mathscr{S}(X)$, then we get an induced hyperstructure on $X$ \textemdash\ $\H(X)$ with $X=B_0$, using the boundary maps in $\mathscr{S}(X)$. (See \newline \citet{B2} for pulling and pushing $\H$-structures.) So
\begin{align*}
    \H(X) = 
    \begin{pmatrix}
    B_n \\
    \vdots \\
    B_0 = X
    \end{pmatrix}
    \rightsquigarrow 
    \begin{pmatrix}
    \mathscr{S}_n \\
    \vdots \\
    \mathscr{S}_0
    \end{pmatrix}
    = \mathscr{S}, \quad \H\text{-states}.
\end{align*}
Alternatively we may have an $\H$-structure where $X=B_n$, with states $\Sigma(X) \subset \mathscr{S}_n$. Then we may assign possible, compatible states to $B_{n-1}$,...
\begin{align*}
    \H(X) =
    \begin{pmatrix}
    B_n=X \\
    B_{n-1} \\
    \vdots \\
    B_0
    \end{pmatrix}
    \begin{matrix}
    \rightsquigarrow \\
    \rightsquigarrow \\
    \vdots \\
    \rightsquigarrow
    \end{matrix} \quad
    \begin{matrix}
    \Sigma_n(X) = \Sigma(X)\subset \mathscr{S}_n \\
    \Sigma_{n-1}(X)\qquad\qquad\qquad \\
    \vdots \qquad\qquad\qquad\\
    \Sigma_0(X) \qquad\qquad\qquad
    \end{matrix}
\end{align*}
Hence via the boundary maps in $\H(X)$, we get the $\Sigma(X)$ states organized into a hyperstructure:
\begin{align*}
    \Sigma(X) = 
    \begin{pmatrix}
    \Sigma_n \\
    \Sigma_{n-1} \\
    \vdots \\
    \Sigma_0 \\
    \end{pmatrix}
\end{align*}
creating new types of lower states ($\H$-states), which may facilitate dynamical changes and actions.

We could also reverse the processes. In the first case take $X=B_n$ and let $\mathscr{S}$ induce an $\H$-str. on $X$ with $B_n$ the top level. In the last case we could take
\begin{align*}
    B_0 = X, \quad \Sigma(X)=\mathscr{S}_0
\end{align*}
and use an $\H$-structure $\H(X)$ given on $X$ to create an $\H$-state structure with $\mathscr{S}_0$ as the bottom level. This follows from using the boundary maps in $\H(X)$, and creates new higher order states on $X$.

This shows that from a given set of states $S=S_0$ we may create higher order states via a hyperstructure $\mathscr{S}(S)$. This is a local to global process. Conversely, $S=S_n$ may be represented at the top level of a hyperstructure $\mathscr{S}(S)$ dissolving the states and creating families of lower level states.

Finally, let us note that forming a hyperstructure often takes resources to be put into the bonds, but breaking them will release resources. This is similar to geometric cutting in hyperstructures of links, see \citet{B1}.

\section{Higher Quantum States \textemdash\ Higher Entanglement}
\label{sec:HQS}

We have already pointed out that hyperstructures are useful in organizing states \textemdash\ including quantum states. In cold gas physics one has detected bound states called Borromean, Brunnian or Efimov states with the property that particles (atoms) bind three by three, but not two by two. Similarly GHZ states for qubits.

These states are having properties corresponding to Borromean and Brunnian links in $\mathbb{R}^3$. In \citet{B1} we introduced, constructed and illustrated ``higher'' order links generalizing Borromean and Brunnian ones.

Thinking of this at the metaphorical level led to the conjecture that there are corresponding higher quantum states which may not yet be detectable by present day technology. The plethora of higher links (see \citet{B1}) makes it natural to suggest the existence of corresponding higher $n$-th order state being organized in the same way as a hyperstructure. We will call these states \textit{$n$-th order entangled states}. GHZ states are of first order and one may along these lines consider higher order GHZ states \textemdash\ for example 2GHZ in analogy with Figure 8. Linking creates the bonds and Obs is the test for circularity and for example the Brunnian property in the higher link hyperstructures.

We conjecture that there is a suitable Topological Quantum Field Theory relating the higher order links and higher order entangled states. This should mean, for example, that the field theory would send $(S^3-nB)$, the complement of the $n$-th order Brunnian link (disjoint union of tori), to an ``$n$-th order' GHZ state $n$-GHZ. We will elaborate on this elsewhere. 

Let us explain higher order entanglement of quantum states in more detail by an example of a composite quantum system of particles: $1,2,3,\ldots,9$ being grouped as in Figure 5:

\begin{figure}[hbt!]
\begin{equation*}
\begin{tikzpicture}
\draw (0,0) circle (1.5);
\draw (4,0) circle (1.5);
\draw (8,0) circle (1.5);
\node[circle, draw=black]at(-0.6,0.4){$1$};
\node[circle, draw=black]at(0.6,0.4){$2$};
\node[circle, draw=black]at(0,-0.6){$3$};
\node[circle, draw=black]at(3.4,0.4){$4$};
\node[circle, draw=black]at(4.6,0.4){$5$};
\node[circle, draw=black]at(4,-0.6){$6$};
\node[circle, draw=black]at(7.4,0.4){$7$};
\node[circle, draw=black]at(8.6,0.4){$8$};
\node[circle, draw=black]at(8,-0.6){$9$};
\node[]at(1.3,-1.3){I};
\node[]at(5.3,-1.3){II};
\node[]at(9.3,-1.3){III};
\draw (-2,-2) rectangle (10,2);
\node[]at(4,2.3){System};
\end{tikzpicture}    
\end{equation*}
Figure 5
\end{figure}

Corresponding state spaces of the particles are:
\begin{align*}
    \H_\text{I} &= \H_1 \otimes \H_2 \otimes \H_3 \\
    \H_\text{II} &= \H_4 \otimes \H_5 \otimes \H_6 \\
    \H_\text{III} &= \H_7 \otimes \H_8 \otimes \H_9 \\
    \H &= \H_\text{System} = \H_\text{I} \otimes \H_\text{II} \otimes \H_\text{III}.
\end{align*}

We consider for simplicity pure states, but the constructions go through for density matrices, mixed states and ensembles as well. 

Pure states are product states:
\begin{align*}
    \text{Pure}(\H_\text{I}) = \{ \alpha_1 \otimes \alpha_2 \otimes \alpha_3 \;\vert\; \alpha_i \in \H_i \},
\end{align*}
other states are called entangled, Ent$(\H_\text{I})$. Similarly for $\H_\text{II}$ and $\H_\text{III}$. We define \newline
\noindent\textit{First level pure states} are of the form
\begin{align*}
    \alpha_1 \otimes \alpha_2 \otimes \ldots \otimes \alpha_9,
\end{align*}
\noindent\textit{First level entangled states} are sums of tensor products
\begin{align*}
    &\sum_{j=1}^k \gamma_j \bigotimes_{i=1}^9 \alpha_i^j, \quad \text{for some } k,\quad \gamma_j \text{ scalars}. \\
    &\\
    \text{Let } A_\text{I} \in \text{Ent}(\H_\text{I}),\; &A_\text{II} \in \text{Ent}(\H_\text{II}),\; A_\text{III} \in \text{Ent}(\H_\text{III})
\end{align*}
Then we suggest: \newline
\noindent \textit{Second level pure states} are of the form
\begin{align*}
    A_\text{I} \otimes A_\text{II} \otimes A_\text{III}.
\end{align*}
\noindent \textit{Second level entangled states} are sums of tensor products of $A_\text{I}$'s, $A_\text{II}$'s and $A_\text{III}$'s. Similarly for more particles and levels. 

Use the composite structure of the system to define higher entanglement. 

The algebraic expressions of sums of tensor products represent the bonds and Obs is the test for pureness and further properties in the hyperstructure these states are organized into. Dissolving a bond means
\begin{align*}
    \partial(n\text{-Entangled state}) = \{ (n-1)\text{-Entangled states} \}.
\end{align*}

Let us conclude this section with some more details. Our general field assignments
\begin{align*}
    \mathcal{F}\colon \H \rightsquigarrow \mathscr{S}
\end{align*}
give the objects in $\H$ higher states in $\mathscr{S}$. In order to guide the general situation let us consider the case when $\mathscr{S}$ is given by higher entangled states. We will construct this hyperstructure
\begin{align*}
    \mathscr{S} = \{\mathscr{S}_1,\mathscr{S}_2,\ldots,\mathscr{S}_n\}
\end{align*}
as follows. We are given a (quantum) system with state space $\mathscr{S}_1$ being a linear convex space with a tensor product. Let us consider composite systems like for example an $n$-particle system with
\begin{align*}
    \mathscr{S}_1 = \H = \H_1 \otimes \H_2 \otimes \cdots \otimes \H_n.
\end{align*}
In $\mathscr{S}_1$ pure states are pure tensorproduct elements, in other contexts corresponding to the extreme points of their convex hull (the Krein-Milman theorem).

We will now sketch the association of a hyperstructure of higher entangled states based on $\mathscr{S}_1$. Let $\mathcal{P}_1 \subset \mathscr{S}_1$ be the pure states, and $P$ stands for finite powersets. We define bonds by
\begin{align*}
    B_1\colon P(\mathcal{P}_1) &\longrightarrow \mathscr{S}_1 \\
    B_1\left(\{q_{i,j}^1\}\right) &= \left( \sum \alpha \; \otimes\right)\left(\{q_{i,j}^1\}\right) \\
    &= \sum_j \alpha_j \otimes_i q_{i,j}^1
\end{align*}
$\alpha_j$'s are scalars.

In the terminology in \cite{B1,B2,B3,B4,B5,B6,B7} we would write
\begin{align*}
    &\qquad S_1 = \{q_{i,j}^1\}, \quad \omega_1 = \text{``property of not being pure'' by Obs}_1. \\
    &\qquad B_1(S_1,\omega_1) = \mathrm{Convex}(S_1) = K_1 \text{ (entangled by $S_1$)}.
\end{align*}
(Corresponding to the convex hull with the extreme points/pure states excluded)

Geometrically:

\begin{center}
\begin{tikzpicture}[scale=0.7]
\draw (0,0) -- (1,2) -- (5,2) -- (6,-1) -- (2,-2) -- (0,0);
\node[] at(3,0){$K_1$};
\node[] at(5.5,2.5){$q_{i,j}^1$};
\node[] at(3,-2.5){Figure 6};
\end{tikzpicture}
\end{center}

We put
\begin{align*}
    \mathrm{Ent}_1 = \mathrm{Im }\ B_1
\end{align*}
by definition the first order entangled elements -- being the bonds and the basic elements of the next level formation. Again $\mathrm{Obs}_2$ checks non-pureness in the sense of not being pure tensor products of elements in $\mathrm{Ent}_1$ (and possibly other chosen properties),
\begin{align*}
    \mathcal{P}_2 = \otimes q_k, \quad q_k \in \mathrm{Ent}_1.
\end{align*}
\begin{align*}
    B_2\colon P(\mathcal{P}_2) &\longrightarrow \mathscr{S}_2 \\
    B_2\left(\{q_{i,j}^2\}\right) &= \left(\sum \alpha \; \otimes \right) \left(\{q_{i,j}^2\}\right)
\end{align*}
Geometrically:

\begin{center}
\begin{tikzpicture}[scale=0.7]
\draw (0,0) -- (1,2) -- (5,2) -- (6,-1) -- (2,-2) -- (0,0);
\draw (1.4,0) -- (2.2,1.25) -- (4.1,1) -- (2.9,-1.2) -- (1.4,0);
\node[] at(4.5,-0.5){$K_1$};
\node[] at(3,0){$K_2$};
\node[] at(5.5,2.5){$q_{i,j}^1$};
\node[] at(4.6,1.2){$q_{i,j}^2$};
\node[] at(3,-2.5){Figure 7};
\end{tikzpicture}
\end{center}

$S_2=\{q_{i,j}^2\}, \; \omega_2 \in \mathrm{Obs}_2$, checking out pure tensorproduct states, but could also include other properties like Brunnianness \cite{B1} or others.
\begin{align*}
    B_2(S_2,\omega_2) = K_2
\end{align*}
In this way we proceed to level $n$ obtaining $B_n$ and $K_n$, the bonding operation being of the type $\left(\sum \alpha \otimes \right)$. But note that the properties in Obs ($\omega_i \in \mathrm{Obs}_i$) may vary and should be chosen with care. We put
\begin{align*}
    \mathrm{Ent}_k = \mathrm{Im}\ B_k
\end{align*}
and call them $k$-th order entangled states.

When we have
\begin{align*}
    B_k\left(\{q_{i,j}^k\}\right) = \left(\sum \alpha \; \otimes \right)\left(\{q_{i,j}^k\} \right) \ni q_{i,j}^{k+1}
\end{align*}
we say that $q_{i,j}^{k+1}$ is entangled by $\{q_{i,j}^k\}$ at level $k$. The states at all levels could be time dependent.

If we put $\mathrm{Ent}_0 = \mathscr{S}_1$ then we have a hyperstructure
\begin{align*}
    \mathrm{Ent} = \{\mathrm{Ent}_0,\mathrm{Ent}_1,\ldots,\mathrm{Ent}_n\}
\end{align*}
of higher order entangled states depending on chosen properties. Field theories of this form may be an interesting object of study leading to systems with new physical properties.

We have here considered a direct definition of entanglement in composite systems in order to illustrate the idea of forming higher entanglements in a simple way. The idea extends to other situations as well.

These higher quantum states with a hyperstructured organization were predicted in \citet{B1}. For example second order entangled states correspond to links of the type in Figure 8 and is a useful metaphor. However, it would be interesting to have these corresponding link states expressed algebraically in the state space.

\begin{figure}[hbt!]
    \begin{center}
        \text{ }\\
        \includegraphics[scale=0.3]{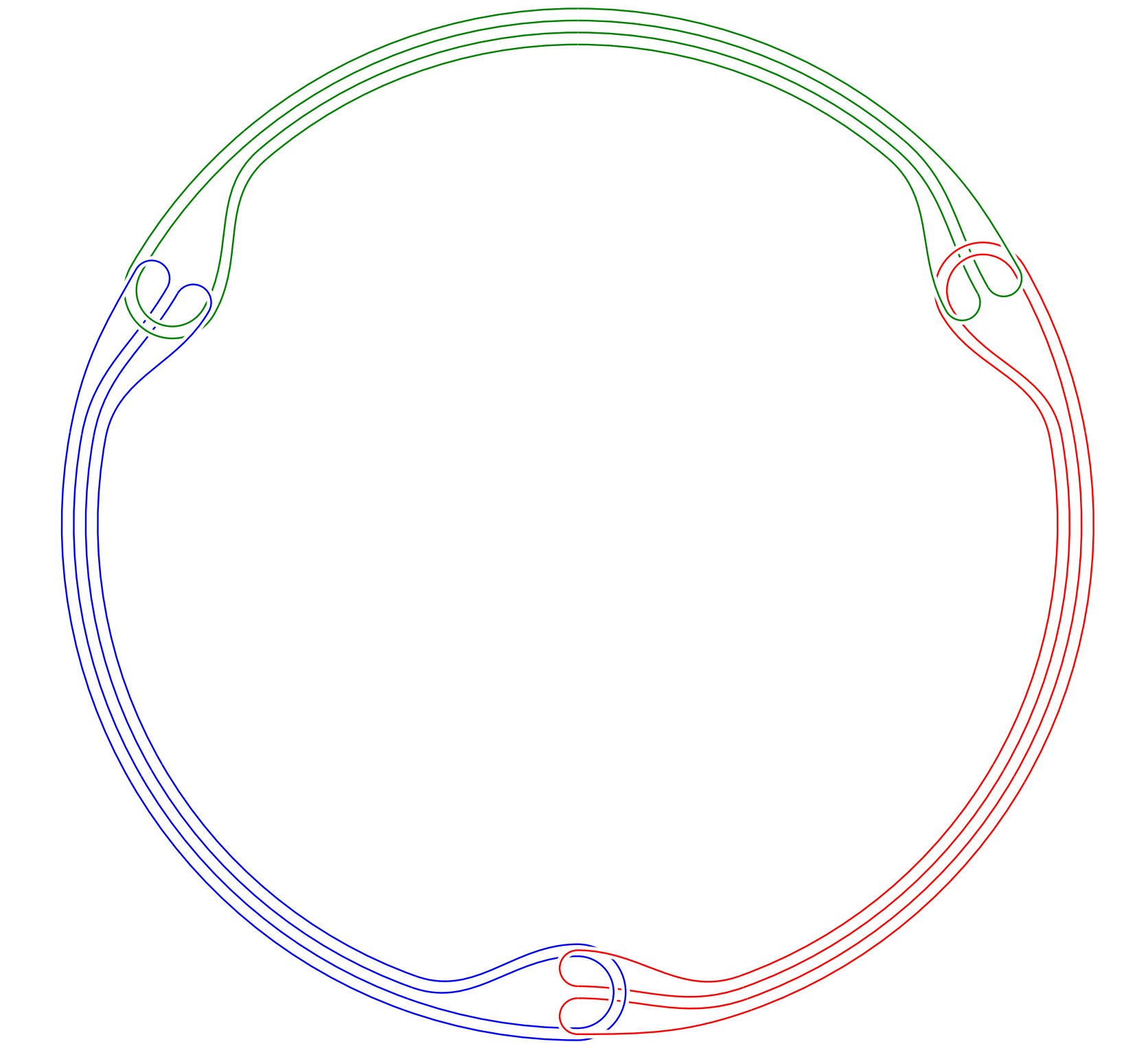} \\ \text{ }\\
        \includegraphics[scale=0.3]{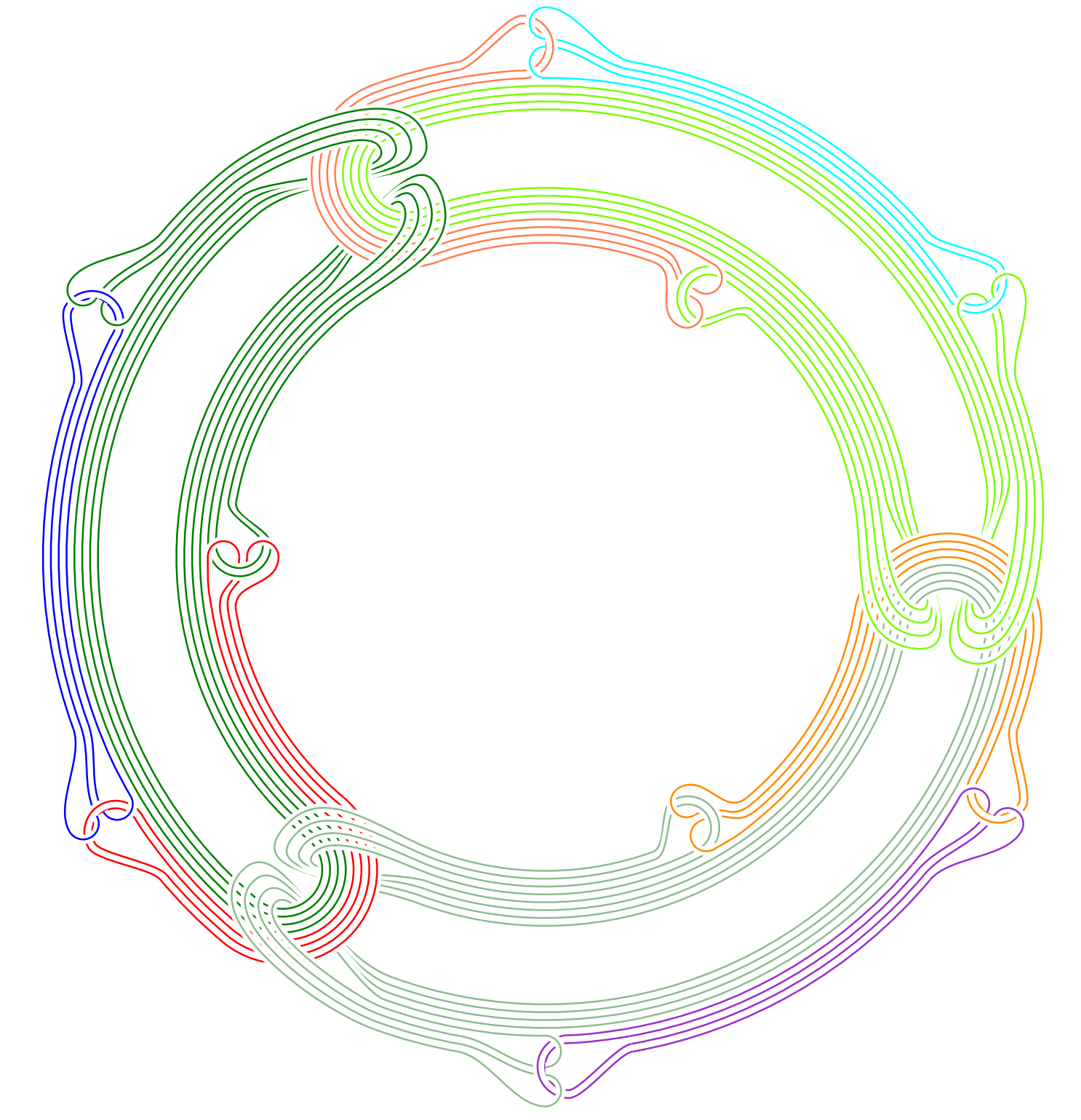}
    \end{center}
    Figure 8
\end{figure}

Other hyperstructures may be induced from this by ``pushing forward'' or ``pulling back'' processes as described in \citet{B2}. Hyperstructures represent a tool to organize states and create higher states.

We may also organize observables into hyperstructures by creating hyperstructures of for example $C^*$-algebras.  

In general hyperstructure may be thought of as ``societies'' of particles, states, objects and structures. We associate fields of states to particles, etc. in the form of a general field theory of hyperstructures. The hyperstructure will then guide the societies of particles that can be formed, like many levels corresponding to ``very social'', highly organized collections of particles, while few levels may reflect ``unsocial'', less organized collections of particles.

We may also use these examples to illustrate ``bond type'' field theories as we have suggested in section 4. \textemdash\ namely $B(\{\H_i\})$.
\begin{align*}
    &\H_1 \text{, may be a hyperstructure of sets of particles.} \\
    &\H_2 \text{, representing the particles by circles in }\mathbb{R}^3 \text{, obtaining higher links.} \\
    &\H_3 \text{, assigning states in a Hilbert space to the particles.}
\end{align*}
Bonds $B(\{H_1,\H_2,\H_3\})$ are the ways in which the hyperstructures interact or relate. These bonds may again interact leading to
\begin{align*}
    \text{Hyp}(\{B_j\{H_i\}\}).
\end{align*}

\newpage

\section{Neural Field Theory}

In neuroscience one studies neurons with individual states \textemdash\ $1$ or $0$ \textemdash\ fire or not. If one collects the firing data for a large colletion $N$ of neurons one is interested in their collective state or behaviour by associating a state space to $N$. Hence we have the situation
\begin{align*}
    &N = \text{ neural population} \\
    &\mathscr{S}(N) = \text{ collective state space} \\
    &N \rightsquigarrow \mathscr{S}(N)
\end{align*}
This has been done for Head Direction Cells \cite{13} and Grid Cells \cite{14}. Let us indicate how this may lead to a Neural Field Theory of the GFT type by further experimental data.

In \cite{15, 16} we have discussed how to organize collections of neurons and their state spaces into hyperstructures and define Neural Field Theories. Let us recall and elaborate some of these ideas. 


In the brain there are lots of modules organized into anatomical and functional layers. The modules at level $n$ we call $n$-modules. Let $N$ be a large collection of neurons (brain cells) \textemdash\ even the whole brain \textemdash\ under consideration. Furthermore, let us assume $N$ is in a natural way \textemdash\ due to anatomy and functionality \textemdash\ divided into subgroups:
\begin{align*}
    N = N_1 \cup N_2 \cup \ldots \cup N_{g_1}
\end{align*}
Du to the level structure of the brain, these subgroups may be further divided into new subgroups
\begin{align*}
    N_i = N(i) = N(i,1)\cup N(i,2) \cup \ldots \cup N(i,g_2).
\end{align*}
This may continue at many stages.

So we may think of $N(i)$ \textit{binding} the $N(i,j)$'s together by \textit{observing} the inclusions using Hyperstructure terminology, see \cite{15,B5,B6}.

\begin{equation*}
\begin{tikzpicture}
\draw (0,0) rectangle (5,2.5);
\draw (1.25,0) -- (1.25,2.5);
\draw (2.5,0) -- (2.5,2.5);
\draw (3.75,0) -- (3.75,2.5);
\draw (2.5,0.5) -- (3.75,0.5);
\draw (2.5,1) -- (3.75,1);
\draw (2.5,1.5) -- (3.75,1.5);
\draw (2.5,2) -- (3.75,2);
\node[]at(4.1,-0.8){$N(i)$};
\draw[->] (3.6,-0.8) arc (-90:-180:0.7);
\node[]at(6.3,0.9){$N(i,j)$};
\draw[->] (6,1.2) arc (20:160:1.5 and 1);
\end{tikzpicture}
\end{equation*}
This iterates for $n$ levels to subgroups parametrized by tuples
\begin{align*}
    \omega &=(i_1,\ldots,i_n), \quad i_k \in I_k \text{ (index set)}. \\
    N(\omega) &= N(i_1,i_2,\ldots,i_n)
\end{align*}
$(i_1,i_2,\ldots,i_n)$ may be thought as ``anatomical'' coordinates, see the Figure.
\begin{align*}
    N = \bigcup_{i_1\in I_1} N(i_1) = \ldots = \bigcup_{\substack{i_1 \in I_1 \\ \vdots \\ i_n \in I_n}} N(i_1,i_2,\ldots,i_n)
\end{align*}
We may also just consider a nesting of $N$ by using inclusions. 

Now we apply the statespace assignment to the collection or assembly
\begin{align*}
    \{N(\omega)\}
\end{align*}
and get an assembly of spaces
\begin{align*}
    \{\mathscr{S}(\omega)\}.
\end{align*}
We assume that all the $\mathscr{S}(\omega)$'s are non-trivial (not the shape of a point), and all the $N(\omega)$ are topological modules. 

\begin{definition}
    A collection of neurons $N$ is called a topological module if its collective state space $\mathscr{S}(N)$ has a ``shape'' different from that of a point.
\end{definition}
\begin{remark}
    In mathematical terms this means that $\mathscr{S}(N)$ is not a contractible space, which is the case both for HD cells and grid cells.
\end{remark}

Furthermore, we assume that the lowest level data point cloud generates a new cloud in the state space and that this continues to higher levels. For grid cells this is like the spike trains determining the path of the rat in the box, see \cite{14}.

The $\mathscr{S}(\omega)$'s are also level structured in the sense that
\begin{align*}
    &\mathscr{S}(i_1,i_2,\ldots,i_{n-1}) \text{ is created from the family} \\
    \{&\mathscr{S}(i_1,i_2,\ldots,i_n)\}_{i_n\in I_n}
\end{align*}
such that the states in $\mathscr{S}(i_1,\ldots,i_{n-1})$ are collective states \textit{binding} together the families of states in $\{\mathscr{S}(i_1,\ldots,i_n)\}_{i_n\in I_n}$ through time \textit{observation}. $I_n$ is some index set. 

The multilevel binding structures
\begin{align*}
    \mathcal{N}=\{N(\omega)\} \text{ and } \mathscr{S}=\{\mathscr{S}(\omega)\}
\end{align*}
are examples of a general framework called \textit{Hyperstructures}, see \cite{15,B5, B6} and references there in for a definition, examples and further information and references.

The assignment
\begin{align*}
    &\mathcal{F}: \mathcal{N} \rightsquigarrow \mathscr{S} \\
    \text{by } &\mathcal{F}(N(\omega)) = \mathscr{S}(\omega)
\end{align*}
is a special case of a general field theory motivated by quantum field theories. Hence, we think that
\begin{align*}
    \mathcal{F}: \mathcal{N} \rightsquigarrow \mathscr{S}
\end{align*}
should be thought of and called a \textit{Neural Field Theory}, see \cite{16}.

Let us write $\mathscr{S} = \{\mathscr{S}_0\}$ as a level structure

\begin{equation*}
\begin{tikzcd}[row sep=tiny]
\mathscr{S}_0 = \mathscr{S}(\emptyset), \text{ top level} \arrow[dd, "\partial"', bend right=35, shift right=17] \\
\\
\mathscr{S}_1 = \mathscr{S}(i_1) \qquad\qquad \text{ } \arrow[uu, "\delta"', bend right=35, shift left=15] \arrow[dd, "\partial"', bend right=35, shift right=17]\\
 \vdots
 \text{ }\\
\mathscr{S}_n = \mathscr{S}(i_1,\ldots,i_n)\quad \text{ } \arrow[uu, "\delta"', bend right=35, shift left=15]
\end{tikzcd}
\end{equation*}
Levelwise coding and decoding via $\partial_i$ and $\delta_i$. 

From the construction of the state spaces we have maps
\begin{align*}
    \partial_i: \mathscr{S}_i \longrightarrow \mathscr{S}_{i+1}
\end{align*}
``dissolving'' the state as a bond. This corresponds to \textit{coding}. Piecing together ``local'' states through levels to ``global'' states via maps $\delta_i$ is what we in \cite{B3,B5, B6} and references there in call a \textit{globalizer}. This corresponds to neural \textit{decoding}.


Let us consider another aspect of neurons. When an animal or basically a brain explores a space there are usually objects in the space, objects of objects, etc. In this context it may be useful to introduce a more abstract setting.
\begin{align*}
    &\text{Space: } X \text{ (set, topological space, manifold,...)} \\
    &\text{Space with objects: } (X, \{V_i\}), \quad V_i \subset X, \\
    &\text{Space with objects of objects ($2$-objects): } (X,\{V_{ij}\}), \\
    &\qquad\qquad\qquad\quad\qquad\qquad\qquad\qquad\qquad\qquad V_{ij}\subset V_i \subset X \\
    & \qquad \vdots \\
    &\text{Space with $n$-objects: } (X,\{V(\omega)\}) \\
    & \omega = (i_1,\ldots,i_n), \quad i_k\in I_k \\
    &\text{ and } V(\omega,i_{n+1}) \subset V(\omega) \\
    & V(\omega) \text{ is a \textit{bond} of } \{V(\omega,i_{n+1})\}
\end{align*}
Hence
\begin{align*}
    \bigcup_\omega V(\omega) \subset X,
\end{align*}
but not necessarily all of $X$. Obs is a property test. $\mathcal{W}=\{V(\omega)\}$ is the multilevel family of objects in the space. If no objects $\mathcal{W}=\{\emptyset\}$. $\mathcal{W}$ gives a hyperstructure on $X$:
\begin{align*}
    \H_\mathcal{W}(X).
\end{align*}
For $\mathcal{W}=\{\text{open subsets}\}$, $\H_\mathcal{W} = \H_\text{Nest}$ in Section 5.

To objects we assign states and properties like in a general field theory
\begin{align*}
    \mathcal{F}: \H_\mathcal{W}(X) \longrightarrow \mathscr{S}
\end{align*}
where $\mathscr{S}$ is a suitable hyperstructure of states and properties. Local properties will be glued together to global properties by a ``globalizer'' mechanism. Of course in specific cases more structure will be added to $X$, $\mathcal{W}$ and $\mathscr{S}$. We think that these structures may be useful in future experiments and the analysis of the data, and representing new types of field theories.

Finally, in many situations algorithms may be organized into hyperstructres. In particular may hyperstructres of neural algorithms in the brain lead to better understanding of higher cognitive functions. In this context hyperstructures may be thought of as a form of deep learning. 

\section{Epilogue}

Hyperstructures are tools for working with structures and situations where several levels come into play. Biological systems -- organisms -- are clearly good examples. Given an object $X$, putting a hyperstructure on it $\H(X)$, makes it in some sense "organismic". On such a hyperstructure it is often useful to assign a field theory of the type we have discussed
\begin{align*}
    &\mathcal{F}\colon \H \rightsquigarrow \mathscr{S} \\
    &\mathscr{S} = \{\mathscr{S}_0,\ldots,\mathscr{S}_n\}
\end{align*}
where $\mathscr{S}_0$ are "global" states and $\mathscr{S}_n$ are "local".

As discussed in \cite{B5,B6} we may introduce a "globalizer" connecting the "local" states in $\mathscr{S}_n$ to a global state in $\mathscr{S}_0$:
\begin{align*}
    \mathrm{Glob}\left(\{z_n\}\right) = z_0
\end{align*}
defined via levels \cite{B2,B5}. 

Often one wants to find an action of the form 
$$\mathcal{A}\times \mathscr{S}\longrightarrow \mathscr{S}$$
changing $z_0$ to $\hat{z}_0$. This may be difficult and take many resources. However, with a field theory as above we may more easily and with less resources find a "local" action taking
\begin{align*}
    \{z_n\} \rightsquigarrow \{\hat{z}_n\}
\end{align*}
such that
\begin{align*}
    \mathrm{Glob}\left\{\hat{z}_n\}\right) = \hat{z}.
\end{align*}
This is "tunnelling through levels".

Often there are barriers to the transition $z_0 \rightsquigarrow \hat{z}_0$. Such barriers may be of very different nature like high energy potentials for example in physical tunnelling phenomena including nuclear fusion. In organisms the barriers may be of cognitive nature. In the superconductors in tunnelling in Josephson junctions the electrons are organized into pairs.

In organisms both cognitive and physical barriers may be overcome by drugs -- acting locally as small perturbations. We suggest that "localizing" the barrier situation by introducing a general field theory may be useful.

\begin{definition}
The described process
\begin{align*}
    z_0 \rightsquigarrow \hat{z}_0
\end{align*}
we call \textit{Abstract Tunnelling (or $\H$-Tunnelling)}.
\end{definition}

We think that this process deserves to be studied in this abstract setting, and that this may be useful in many physical situations, decoding of local neural brain activity, see \cite{14,16}, and switching from one attractor to another.

\begin{sloppypar}$\H$-structures are really structures where object composition takes place through bonds (synthesis, fusion). Composed objects require bindings but sometimes less when suitably organized (like an $\H$-structure), and bindings may be released in a desired way. This applies to energy, objects, properties (physical, mathematical, ...). In a hyperstructure one may also do structured decomposition through dissolving bonds (analysis, fission) leading to positive or negative releases as for composition depending on the starting objects.\end{sloppypar}

When we want to act on an object $X$ (a set or a collection) we may do this through the use of field theories of hyperstructures. Let us start with $X$ and we choose a hyperstructure
\begin{align*}
    \H = \{H_0,H_1, \ldots, H_n\}
\end{align*}
suitable for the situation we have in mind, The $H$'s are the underlying basic level sets. We choose a representative map
\begin{align*}
    \varphi \colon X = X_0 \longrightarrow H_0 \quad \text{(lowest level)}
\end{align*}
If $S_0 \subset X_0$, $\varphi(S_0) \subset H_0$ and we put $\widehat{B}_0(S_0) = B_0(\varphi(S_0))$ and similarly for higher bonds giving $\widehat{B}_i$. The $B$'s are bonds in $\H$ and the $\widehat{B}_i$'s are the new bonds in the induced hyperstructure $\H(X)$ on $X$. This process is described in \cite{B2}.

We pass from the original object $X$ to an associated hyperstructure $\H(X)$. The next step is then to assign fields ("states") that will enable us to act on $\H(X)$ and $X$ in a desired way
\begin{align*}
    \mathcal{F}\colon \H(X) \rightsquigarrow \mathscr{S}
\end{align*}
where $\mathscr{S}$ is a hyperstructure representing the fields. In the framework of $\H$-fields we also get $\H$-forces and $\H$-charges.

\end{document}